\documentclass[draft]{article}
\usepackage{amssymb,amsmath}
\usepackage{graphics}
\input xy
\xyoption{all}

\begin{document}
\title{Complex of twistor operators in symplectic spin geometry}
\author{Svatopluk Kr\'ysl \footnote{{\it E-mail address}: krysl@karlin.mff.cuni.cz}\\ {\it \small  Charles University of Prague, Sokolovsk\'a 83, Praha,  Czech Republic}.
\thanks{I thank to V\'{i}t Tu\v{c}ek and Andreas \v{C}ap for discussions and comments.
The author of this article was supported by the grant
GA\v{C}R 306-33/52117 of the Grant Agency of Czech Republic. The work is a part of the 
research project MSM 0021620839 financed by M\v{S}MT \v{C}R.}}

\maketitle \noindent
\centerline{\large\bf Abstract} For  a symplectic manifold admitting a metaplectic structure (a symplectic analogue of the Riemannian spin structure), we construct a sequence consisting of differential operators using a symplectic torsion-free affine connection. All but one of these operators are of first order.  The first order ones are symplectic analogues of the twistor operators known from Riemannian spin geometry. We prove that under the condition the symplectic Weyl curvature tensor field of the symplectic connection vanishes, the mentioned sequence forms a complex.
This gives rise to a new complex for the so called Ricci type symplectic manifolds, which admit a metaplectic structure.

{\it Math. Subj. Class.:} 53C07, 53D05, 58J10.

{\it Keywords:} Fedosov manifolds, metaplectic structures, symplectic spinors, Kostant spinors, Segal-Shale-Weil representation, complexes of differential operators.

\section{Introduction}

  In the paper, we shall introduce a sequence of   differential operators acting on symplectic spinor valued exterior differential forms over a symplectic manifold $(M,\omega)$ admitting the so called metaplectic structure.  To define these operators, we make use of a symplectic torsion-free affine connection $\nabla$ on $(M,\omega).$
Under certain condition on the curvature of the connection $\nabla,$ described bellow, we prove that the mentioned sequence forms a complex.

  Let us say a few words about the metaplectic structure.
The symplectic group
$Sp(2l,\mathbb{R})$ admits a non-trivial two-fold covering, the so called metaplectic group, which we shall  denote  by $Mp(2l,\mathbb{R}).$ Let $\mathfrak{g}$ be the Lie algebra of $Mp(2l,\mathbb{R}).$
A~metaplectic structure on a symplectic manifold $(M^{2l},\omega)$ is
a notion parallel to a spin structure on a Riemannian manifold.  In particular,
one of its part is a principal $Mp(2l,\mathbb{R})$ bundle $(q: \mathcal{Q} \to M,Mp(2l,\mathbb{R}))$.

For a symplectic manifold admitting a metaplectic structure, one can construct the so called symplectic spinor bundle $\mathcal{S} \to M,$  introduced  by Bertram Kostant in 1974.  The symplectic spinor bundle $\mathcal{S}$ is the vector bundle associated to the metaplectic structure $(q: \mathcal{Q}\to M,Mp(2l,\mathbb{R}))$
on $M$ via the so called  Segal-Shale-Weil representation
of the metaplectic group $Mp(2l,\mathbb{R})$. 
See Kostant \cite{Kostant2} for details. 
 
  The Segal-Shale-Weil representation is an infinite dimensional unitary representation of
the metaplectic group $Mp(2l, \mathbb{R})$ on the space of all complex valued square
Lebesgue integrable functions ${\bf L^{2}}(\mathbb{R}^{l}).$ Because of the infinite dimension, the  Segal-Shale-Weil representation is not so easy to handle. It is known, see, e.g., Kashiwara, Vergne \cite{KV}, that  
the $\mathfrak{g}^{\mathbb{C}}$-module structure of the underlying Harish-Chandra module of
this representation is equivalent to the space
$\mathbb{C}[x^1,\ldots, x^l]$ of polynomials in $l$ variables, on which
the Lie algebra $\mathfrak{g}^{\mathbb{C}}\simeq \mathfrak{sp}(2l,\mathbb{C})$ acts 
via the so called Chevalley homomorphism,\footnote{The Chevalley homomorphism
is a Lie algebra monomorhism of the complex symplectic Lie algebra $\mathfrak{sp}(2l,\mathbb{C})$ into the Lie algebra of the
associative algebra of polynomial coefficients differential operators acting on
$\mathbb{C}[x^1,\ldots, x^l].$} see Britten, Hooper, Lemire
\cite{BHL}. Thus, the infinitesimal structure of the Segal-Shale-Weil
representation can be viewed as the complexified {\it symmetric}
 algebra $(\bigoplus_{i=0}^{\infty}\odot^i
\mathbb{R}^l)\otimes_{\mathbb{R}}{\mathbb{C}} \simeq \mathbb{C}[x^1,\ldots, x^l]$ of the Lagrangian
subspace $(\mathbb{R}^l,0)$ of the canonical symplectic vector space
$\mathbb{R}^{2l}\simeq (\mathbb{R}^l,0)\oplus (0,\mathbb{R}^l).$
This shows that the situation is completely parallel to
the complex orthogonal case, where the spinor representation can be
realized as the {\it exterior} algebra of a maximal isotropic
subspace.  An interested reader is
referred to Weil \cite{Weil}, Kashiwara, Vergne \cite{KV} and also
to Britten, Hooper, Lemire \cite{BHL} for details.  For some technical reasons, we shall be using the so called minimal globalization of
the underlying  Harish-Chandra module of the  Segal-Shale-Weil representation, which we will call {\it metaplectic representation} and  denote it by ${\bf S}.$ The elements of $\bf S$ will be called symplectic spinors. 

Now, let us consider a symplectic manifold $(M,\omega)$ together with 
 a symplectic torsion-free affine connection $\nabla$ on it. Such connections are usually called Fedosov connections.
 Because the Fedosov connection is not unique for a choice of $(M,\omega)$ (in the contrary to Riemannian geometry), it seems natural to add the connection to the studied symplectic structure and
 investigate  the triples $(M,\omega, \nabla)$ consisting of a symplectic manifold $(M,\omega)$ and a Fedosov connection $\nabla.$ Such triples are usually called Fedosov manifolds and they were used in the deformation quantization. See, e.g., Fedosov \cite{Fedosov}. Let us recall that in Vaisman \cite{Vaisman}, the space of the so called symplectic curvature tensors was decomposed wr. to $Sp(2l,\mathbb{R}).$ For $l=1,$ the module of symplectic curvature tensors is irreducible, while for $l\geq 2,$ it decomposes into two irreducible submodules. These modules are usually called symplectic Ricci and symplectic Weyl modules, respectively. This decomposition translates to differential geometry level giving rise to the symplectic Ricci and symplectic Weyl curvature tensor fields, which add up to the curvature tensor field of $\nabla.$ See Vaisman \cite{Vaisman} and also Gelfand, Retakh, Shubin \cite{GSR} for a comprehensive treatment on Fedosov manifolds.

Now, let us suppose that a Fedosov manifold $(M,\omega,\nabla)$ admits a metaplectic structure $(q:\mathcal{Q}\to M^{2l}, Mp(2l,\mathbb{R})).$  Let $\mathcal{S} \to M$ be the symplectic spinor bundle  associated to $(q:\tilde{Q} \to M,Mp(2l,\mathbb{R}))$ and let us consider the space $\Omega^{\bullet}(M,\mathcal{S})$ of exterior differential forms with values in $\mathcal{S},$ i.e.,  $\Omega^{\bullet}(M,\mathcal{S}):=\Gamma(M,\mathcal{Q}\times_{\rho}(\bigwedge^{\bullet}(\mathbb{R}^{2l})^*\otimes {\bf S})),$ where $\rho$  is the obvious tensor product representation of $Mp(2l,\mathbb{R})$ on $\bigwedge^{\bullet}(\mathbb{R}^{2l})^*\otimes {\bf S}.$ 
In Kr\'ysl \cite{KryslSVF}, the $Mp(2l,\mathbb{R})$-module $\bigwedge^{\bullet}(\mathbb{R}^{2l})^*\otimes {\bf S}$ was decomposed into irreducible submodules. The elements of $\bigwedge^{\bullet}(\mathbb{R}^{2l})^* \otimes {\bf S}$ are specific examples of the so called higher symplectic spinors. For $i=0,\ldots, 2l,$ let us denote the so called Cartan component of the tensor product $\bigwedge^{i}(\mathbb{R}^{2l})^*\otimes {\bf S}$ by ${\bf E}^{i m_i}.$  (For $i=0,\ldots, 2l,$ the numbers $m_i$ will be specified in the text.)
For $i=0,\ldots 2l-1,$ we introduce an operator $T_i$ acting between the sections of the vector bundle $\mathcal{E}^{im_i}$ associated to ${\bf E}^{i m_i}$ and the sections of the vector bundle $\mathcal{E}^{i+1,m_{i+1}}$associated to ${\bf E}^{i+1,m_{i+1}}.$ In a parallel to the Riemannian case, we shall call these operators symplectic twistor operators.
These operators are first order differential operators and they are defined using the symplectic torsion-free affine connection $\nabla$ as follows. First, the connection $\nabla$ induces a covariant derivative $\nabla^S$ on the bundle $\mathcal{S} \to M$ in the usual way. Second, the covariant derivative $\nabla^S$ determines the associated exterior covariant derivative, which we denote  by $d^{\nabla^S}.$ For $i=0,\ldots, 2l-1$, we define the symplectic twistor operator  $T_i$ as the restriction of $d^{\nabla^S}$ to $\Gamma(M,\mathcal{E}^{im_i})$ composed with the projection to $\Gamma(M,\mathcal{E}^{i+1,m_{i+1}}).$

Because we would like to derive a condition under which $T_{i+1}T_i=0,$ $i=0,\ldots, 2l-1,$ we should focus
our attention to the curvature tensor $R^{\Omega^{\bullet}(M,\mathcal{S})}:=d^{\nabla^S}d^{\nabla^S}$ of $d^{\nabla^S}$ acting on the space $\Omega^{\bullet}(M,\mathcal{S}).$ The curvature $R^{\Omega^{\bullet}(M,\mathcal{S})}$ depends only on the curvature of the symplectic connection $\nabla,$ which consists of the symplectic Ricci and symplectic Weyl curvature tensor fields as we have already mentioned.  In the paper, we will analyze the action of the symplectic Ricci curvature tensor field  on symplectic spinor valued exterior differential forms and especially on $\Gamma(M,\mathcal{E}^{i m_i}),$ $i=0,\ldots, 2l-2.$ We shall prove that the symplectic Ricci curvature tensor field when restricted to $\Gamma(M,\mathcal{E}^{im_i})$ maps this submodule into at most three $Mp(2l,\mathbb{R})$-submodules sitting 
in symplectic spinor valued forms of degree $i+2,$ $i=0,\ldots, 2l-2.$
These  submodules  will be explicitly described.  This will help us to prove that $T_{i+1}T_{i}=0$ $(i=0,\ldots, l-2)$ and 
$T_{i+1}T_i=0$ $(i=l,\ldots, 2l-2)$ assuming the symplectic Weyl curvature tensor field vanishes.   In this way, we will obtain two complexes.
Unfortunately,  one can not expect $T_lT_{l-1}=0$ in general.
This will influence the way, how we construct one complex of the two complexes introduced above.
Let us notice that similar complex was investigated in Severa \cite{Severa} in the case of spheres equipped with the  conformal structure of their round metrics.

The reader interested in applications of the symplectic spinor fields in theoretical physics is
referred to Green, Hull \cite{GH}, where the symplectic spinors are
used in the context of 10 dimensional super string theory. In Reuter
\cite{Reuter}, symplectic spinors are used in the theory of the so called
Dirac-K\"{a}hler fields.

In the second section, some   basic facts on
 the metaplectic representation and higher symplectic spinors are recalled. In this section, we also introduce several mappings acting on the graded space $\bigwedge^{\bullet}(\mathbb{R}^{2l})^*\otimes {\bf S},$ derive some \linebreak (super-)~commutation relations between them and determine a superset of the image of two of them, which are components of an infinitesimal version of the symplectic Ricci curvature tensor field.  
 In the section 3, basic properties of  torsion-free symplectic connections and their curvature tensor field are mentioned and the metaplectic structure is introduced. In the subsection 3.1., the theorem on the complex consisting of the symplectic twistor operators is presented and proved.

\section{Metaplectic representation, higher symplectic spinors and basic notation}
 
 To fix a notation, let us recall some notions from   symplectic linear algebra.
Let us consider a real symplectic vector space
$(\mathbb{V},\omega)$ of dimension $2l,$ i.e., $\mathbb{V}$ is a
$2l$ dimensional real vector space and $\omega$ is a
non-degenerate antisymmetric bilinear form on $\mathbb{V}.$ Let us
choose two Lagrangian subspaces\footnote{maximal isotropic wr. to $\omega$} $\mathbb{L}, \mathbb{L}' \subseteq
\mathbb{V}$ such that $\mathbb{L}\oplus \mathbb{L}'=\mathbb{V}.$ It follows that $\mbox{dim}(\mathbb{L})=\mbox{dim}(\mathbb{L}')=l.$
Throughout this article, we shall use a symplectic basis
$\{e_i\}_{i=1}^{2l}$ of $\mathbb{V}$ chosen in such a way   that
$\{e_i\}_{i=1}^l$ and $\{e_i\}_{i=l+1}^{2l}$ are respective bases of
$\mathbb{L}$ and $\mathbb{L}'.$ Because the definition of  a symplectic
basis is not unique, let us fix   one which   shall be  used in this
text. A basis $\{e_i\}_{i=1}^{2l}$ of $\mathbb{V}$ is called
symplectic basis of $(\mathbb{V},\omega)$ if
$\omega_{ij}:=\omega(e_i,e_j)$ satisfies $\omega_{ij}=1$ if and
only if $i\leq l$ and $j=i+l;$ $\omega_{ij}=-1$ if and only if $i>l$
and $j=i-l$ and finally, $\omega_{ij}=0$ in other cases. Let
$\{\epsilon^i\}_{i=1}^{2l}$ be the basis of $\mathbb{V}^*$ dual to
the basis $\{e_i\}_{i=1}^{2l}.$  For $i,j=1,\ldots,
2l,$ we define $\omega^{ij}$ by
$\sum_{k=1}^{2l}\omega_{ik}\omega^{jk}=\delta_i^j,$ for $i,j=1,\ldots,
2l.$ Notice that not only $\omega_{ij}=-\omega_{ji},$ but also
$\omega^{ij}=-\omega^{ji},$ $i,j =1, \ldots, 2l.$  

 As in the orthogonal case, we would like to rise and lower indices.
Because the symplectic form $\omega$ is antisymmetric, we should be more careful in this case. 
For coordinates ${K_{ab\ldots c\ldots d}}^{rs \ldots t \ldots u}$ of a tensor $K$ over $\mathbb{V},$  we denote
the expression $\omega^{ic}{K_{ab\ldots c \ldots d}}^{rs \ldots t}$ by 
${{{K_{ab \ldots}}^{i}}_{\ldots d}}^{rs \ldots t}$ and 
${K_{ab\ldots c}}^{rs \ldots t \ldots u}\omega_{ti}$ by ${{{K_{ab \ldots c}}^{rs\ldots}}_{i}}^{\ldots u}$ and similarly for other types of tensors and also in the geometric setting when we will be considering tensor fields over a symplectic manifold $(M,\omega)$.   

 Let us denote the symplectic group of $(\mathbb{V},\omega)$  by
$G,$ i.e., $G :=Sp(\mathbb{V},\omega)\simeq Sp(2l,\mathbb{R}).$
Because the maximal compact subgroup $K$ of $G$ is isomorphic to the
unitary group $K \simeq U(l)$ which is of homotopy type
$\mathbb{Z},$ there exists a nontrivial two-fold covering
$\tilde{G}$ of $G.$ See, e.g., Habermann, Habermann \cite{HH} for details. This two-fold covering is called metaplectic
group  of $(\mathbb{V},\omega)$ and it is denoted by
$Mp(\mathbb{V},\omega)$.  Let us remark that $Mp(\mathbb{V},\omega)$ is reductive in the sense of Vogan \cite{Vogan}.
In the considered case, we have
$\tilde{G}\simeq Mp(2l,\mathbb{R}).$ For a later use, let us reserve
the symbol $\lambda$ for the mentioned covering. Thus $\lambda:
\tilde{G} \to G$  is a fixed member  of the isomorphism class of all
nontrivial  $2:1$ covering homomorphisms of $G$.
 Because $\lambda:\tilde{G}\to G$
is a homomorphism of Lie groups and $G$ is a subgroup of the general
linear group $GL(\mathbb{V})$ of $\mathbb{V},$ the mapping $\lambda$ is
also a representation of the metaplectic group $\tilde{G}$ on the
vector space $\mathbb{V}.$ Let us define $\tilde{K}:=\lambda^{-1}(K).$ Obviously, $\tilde{K}$ is a maximal compact subgroup of $\tilde{G}.$
 Further, one can easily see that $\tilde{K}\simeq \widetilde{U(l)}:=\{(g,z)\in U(l)\times \mathbb{C}^{\times}| \mbox{det}(g)=z^2\}$ and thus in particular, $\tilde{K}$ is connected. The Lie algebra $\tilde{\mathfrak{g}}$ of $\tilde{G}$ is isomorphic to the Lie algebra $\mathfrak{g}$ of $G$ and we will identify them. One has $\mathfrak{g}=\mathfrak{sp}(\mathbb{V},\omega)\simeq \mathfrak{sp}(2l,\mathbb{R}).$

 Now let us recall some notions from representation theory which we shall need in this paper. From the point of view of this article, these notions are rather of a technical character.  
 Let $\mathcal{R}(\tilde{G})$ be the category the object of which are complete, locally convex, Hausdorff topological spaces with a  continuous linear $\tilde{G}$-action, such that the resulting representation is admissible and of finite length; the morphisms are continuous $\tilde{G}$-equivariant linear maps between the objects. Let $\mathcal{HC}(\mathfrak{g},\tilde{K})$ be the category of Harish-Chandra $(\mathfrak{g},\tilde{K})$-modules and let us consider the forgetful Harish-Chandra functor $HC:\mathcal{R}(\tilde{G})\to \mathcal{HC}(\mathfrak{g},\tilde{K}).$
It is well known that there exists an adjoint functor $mg: \mathcal{HC}(\mathfrak{g},\tilde{K})\to \mathcal{R}(\tilde{G})$ to the Harish-Chandra functor $HC$. This functor is usually called the minimal globalization  functor and its existence is a deep result in representation theory. For details and for the existence of the minimal globalization functor $mg,$ see Kashiwara, Schmid \cite{KS} or Vogan \cite{Vogan}. 
 
 From now on, we shall restrict ourselves to the case $l\geq 2$
 not alway mentioning it explicitly. The case $l=1$ should be handled separately (though analogously) because
  the shape of the root system of $\mathfrak{sp}(2,\mathbb{R})\simeq \mathfrak{sl}(2,\mathbb{R})$ is different from that
  one of
  of the root system of $\mathfrak{sp}(2l,\mathbb{R})$ for $l\geq 2.$
   As usual, we shall denote the
complexification of $\mathfrak{g}$ by $\mathfrak{g}^{\mathbb{C}}.$
Obviously, $\mathfrak{g}^{\mathbb{C}}\simeq
\mathfrak{sp}(2l,\mathbb{C}).$ 
 
Further, for any Lie group $G$ and a principal $G$-bundle $(p:\mathcal{P} \to M,G)$ over a manifold $M,$ we shall denote the vector bundle associated to this
principal bundle via a representation $\sigma: G \to
\hbox{Aut}(\bf W)$ of $G$ on ${\bf W}$ by $\mathcal{W},$ i.e.,
$\mathcal{W}=\mathcal{G}\times_{\sigma} {\bf W}.$
Let us also mention that we shall often use the Einstein summation convention for repeated indices (lower and upper) without mentioning it explicitly.

\subsection{Metaplectic representation and symplectic spinors}

There exists a distinguished infinite  dimensional unitary representation of
the metaplectic group $\tilde{G}$ which does not descend to a
representation of the symplectic group $G.$ This representation,
called {\it Segal-Shale-Weil},\footnote{The names oscillator or
metaplectic representation are also used in the literature. We shall
use the name Segal-Shale-Weil  in this text, and reserve the name
metaplectic for certain representation arising from the
Segal-Shale-Weil one.} plays an important role in geometric
quantization of Hamiltonian mechanics, see, e.g., Woodhouse
\cite{Wood}. We shall not give a
definition of this representation here and  refer the interested
reader to Weil \cite{Weil} or  Habermann, Habermann \cite{HH}.
 
 The Segal-Shale-Weil representation, which we shall denote by $U,$ is a complex infinite dimensional unitary representation of
$\tilde{G}$ on the space of complex valued square Lebesgue
integrable functions defined on the Lagrangian subspace
$\mathbb{L},$ i.e.,
$$U: \tilde{G} \to
\mathcal{U}({\bf L^2}(\mathbb{L})),$$ where $\mathcal{U}({\bf W})$
denotes the group of unitary operators on a Hilbert space ${\bf W}.$
In order to be precise, let us refer to the space ${\bf
L^2}(\mathbb{L})$ as to the Segal-Shale-Weil module. It is known  that the Segal-Shale-Weil module belongs to the category $\mathcal{R}(\tilde{G}).$ (See   Kashiwara, Vergne \cite{KV} for details and Segal-Shale-Weil representation in general.)
It is easy to see
that the Segal-Shale-Weil representation  splits into two
irreducible $Mp(2l,\mathbb{R})$-submodules ${\bf L^{2}}(\mathbb{L})\simeq {\bf
L^{2}}(\mathbb{L})_+\oplus {\bf L^{2}}(\mathbb{L})_-.$ The first
module consists of even and the second one of  odd complex valued square
Lebesgue integrable  functions on the Lagrangian subspace
$\mathbb{L}.$ Let us remark that a typical construction of the
Segal-Shale-Weil representation is based on the so called
Schr\"{o}dinger representation of the Heisenberg group of
$(\mathbb{V}=\mathbb{L}\oplus\mathbb{L}',\omega)$ and a use of the
Stone-von Neumann theorem.
 
 For technical reasons, we shall need the minimal
globalization of the underlying Harish-Chandra $(\mathfrak{g},\tilde{K})$-module $HC({\bf L^2}(\mathbb{L}))$  of the introduced Segal-Shale-Weil module.  We
shall call this minimal globalization  {\it metaplectic representation} and
denote it by $meta,$ i.e.,
$$meta: \tilde{G} \to \hbox{Aut}(mg(HC({\bf L^2}(\mathbb{L})))),$$ where $mg$ is the minimal globalization functor (see this section and the references therein). For our convenience, let us denote the module
$mg(HC({\bf L^2}(\mathbb{L})))$ by ${\bf S}.$ Similarly we define $\bf S_+$ and $\bf S_-$
to be the minimal globalizations of the underlying Harish-Chandra $(\mathfrak{g},\tilde{K})$-modules of the modules
${\bf L^2}(\mathbb{L})_+$ and ${\bf L^{2}}(\mathbb{L})_-.$
Accordingly to ${\bf L^{2}}(\mathbb{L})\simeq {\bf L^{2}}(\mathbb{L})_+\oplus
{\bf L^{2}}(\mathbb{L})_-,$  we have $\bf S \simeq {\bf S_+} \oplus {\bf S_-}.$
We shall call the $Mp(\mathbb{V},\omega)$-module   $\bf S$
the  symplectic spinor module  and its elements  {\it symplectic spinors}.  For
the name "spinor", see Kostant \cite{Kostant2} or the Introduction.

Further notion related to the symplectic vector space
$(\mathbb{V}=\mathbb{L}\oplus \mathbb{L}',\omega)$ is the so called  symplectic Clifford
multiplication of elements of ${\bf S}$ by vectors from $\mathbb{V}.$ For $i=1,\ldots, l$ and a symplectic spinor $f\in {\bf S},$ we define 
\begin{eqnarray*}
(e_i.f)(x)&:=&\imath x^i f(x) \mbox{ and}\\
(e_{i+l}.f)(x)&:=&\frac{\partial f}{\partial x^{i}}(x),
\end{eqnarray*} where $x=\sum_{i=1}^{l}x^i e_i \in \mathbb{L}$ and $\imath=\sqrt{-1}$ denotes the imaginary unit. 
Extending
this multiplication $\mathbb{R}$-linearly, we get the mentioned
symplectic Clifford multiplication.
Let us mention that the multiplication and the differentiation make sense for any $f\in {\bf S}$ because of the "analytic" interpretation of the minimal globalization. (See Vogan \cite{Vogan} for details.) Let us remark that in the physical literature, the symplectic Clifford multiplication is usually called  Schr\"odinger quantization prescription.
 
The following lemma is an easy consequence of the definition of the symplectic Clifford multiplication.

{\bf Lemma 1:} For $v,w \in \mathbb{V}$ and $s \in {\bf S},$ we have
$$v.(w.s)-w.(v.s)=-\imath \omega(v,w)s.$$

{\it Proof.} See Habermann, Habermann \cite{HH}, pp. 11. $\Box$

 Sometimes, we shall write $v.w.s$ instead of $v.(w.s)$ for $v,w \in \mathbb{V}$ and a symplectic spinor $s\in {\bf S}$ and similarly for higher number of multiplying elements. Further instead of $e_i.e_j.s,$ we shall write $e_{ij}.s$ simply 
and similarly for expressions with higher number of multiplying elements, e.g., $e_{ijk}.s$ abbreviates $e_i.e_j.e_k.s.$
 
\subsection{Higher symplectic spinors}

In this subsection, we shall present a result on a decomposition
of the tensor product of the metaplectic representation  $meta: \tilde{G} \to \mbox{Aut}({\bf S})$ with the  wedge
power of the representation $\lambda^*: \tilde{G} \to
GL(\mathbb{V}^*)$ of $\tilde{G}$ (dual to the representation
$\lambda$) into irreducible summands.  
Let us reserve the symbol $\rho$ for the
mentioned  tensor product representation of $\tilde{G}$, i.e.,\begin{eqnarray*}
&&\rho: \tilde{G} \to \hbox{Aut}(\bigwedge ^{\bullet}\mathbb{V}^*\otimes {\bf S})\\
&&\rho(g)(\alpha\otimes s):=\lambda(g)^{*\wedge r}\alpha\otimes
meta(g)s
\end{eqnarray*}
 for $r = 0,\ldots, 2l,$ $g\in \tilde{G},$   $\alpha\in \bigwedge^r\mathbb{V}^*,$  $s\in \bf{S}$ and extend it linearly.  For definiteness, let us equip the tensor product
$\bigwedge^{\bullet}\mathbb{V}^*\otimes {\bf S}$ with the so called
Grothendieck tensor product topology. See Vogan \cite{Vogan} and Treves \cite{Treves} for
details on this topological structure. In a parallel to the Riemannian case, we shall call the elements of $\bigwedge^{\bullet}\mathbb{V}^* \otimes {\bf S}$  higher symplectic spinors.

Let us introduce the following subsets of the set of pairs of non-negative integers. 
We define
\begin{eqnarray*} 
&&\Xi:=\{(i,j)\in \mathbb{N}_0 \times \mathbb{N}_0 | i = 0, \ldots, l;j=0,\ldots, i \} \cup \\
&& \qquad \mbox{  } \cup \{(i,j)\in \mathbb{N}_0\times \mathbb{N}_0 |i=l+1,\ldots, 2l, j=0,\ldots, 2l-i\},\\ 
&&\Xi_+:=\Xi - \{(i,i)|i=0, \ldots, l\} \, \mbox{ and}\\
&&\Xi_-:=\Xi - \{(i,2l-i)|i=l,\ldots, 2l\}.
\end{eqnarray*}

For each $(i,j) \in \Xi,$ a $\mathfrak{g}^{\mathbb{C}}$-module $\mathbb{E}^{ij}_{\pm}$ was introduced in Kr\'ysl \cite{KryslSVF}. These modules are irreducible infinite dimensional highest modules over $\mathfrak{sp}(\mathbb{V},\omega)^{\mathbb{C}}$ and they are described via their highest weights in the mentioned article.
In the next theorem,  the module of symplectic spinor valued exterior forms $\bigwedge^{\bullet}\mathbb{V}^*\otimes {\bf S}$ is decomposed into irreducible submodules.

{\bf Theorem 2:} For $l\geq 2,$ the following decomposition into irreducible \linebreak $Mp(\mathbb{V},\omega)$-submodules
$$\bigwedge^{i}\mathbb{V}^*\otimes {\bf S}_{\pm} \simeq \bigoplus_{j, (i,j) \in \Xi} {\bf E}^{ij}_{\pm}, \quad i=0,\ldots, 2l, \, \mbox{  holds.}$$ 
The modules ${\bf E}^{ij}_{\pm}$ are determined, as objects in the category $\mathcal{R}(\tilde{G}),$ by the fact that first they are submodules of the corresponding tensor product and second the $\mathfrak{g}^{\mathbb{C}}$-structure of $HC({\bf E}^{ij}_{\pm})$ is isomorphic to $\mathbb{E}^{ij}_{\pm}.$

{\it Proof.} See Kr\'ysl \cite{KryslSVF} or Kr\'ysl \cite{KryslJRT}. $\Box$

In the Figure 1, the decomposition in the case $l=3$ is displayed.
In the $i^{th}$ column of the Figure 1, when counted from zero, the summands of $\bigwedge^{i}\mathbb{V}^*\otimes {\bf S},$ $i=0,\ldots,6,$ are written. The meaning of the arrows at the figure will be explained later. 

{\bf Remark:} Let us mention that for any $(i,j), (i,k) \in \Xi,$ $j\neq k,$ we have $\mathbb{E}^{ij}_{\pm} \not\simeq \mathbb{E}^{ik}_{\pm}$ (as $\mathfrak{g}^{\mathbb{C}}$-modules) for all combinations of $\pm$ on the left hand as well as on the right hand side. Using this fact, we have that for $i =0,\ldots, 2l$ the $\tilde{G}$-modules $\bigwedge^i\mathbb{V}^*\otimes {\bf S}_{\pm}$ are multiplicity free. Moreover for $(i,j), (k,j) \in \Xi$, we have $\mathbb{E}^{ij}_{\pm} \simeq \mathbb{E}^{kj}_{\mp}.$ These facts will be crucial in the paper.

For our convenience, let us set ${\bf E}^{ij}_{\pm}:=\{0\}$ for $(i,j) \in \mathbb{Z}\times \mathbb{Z} - \Xi$ and ${\bf E}^{ij}:={\bf E}^{ij}_+\oplus {\bf E}^{ij}_-.$ 

$$\xymatrix{{\bf E}^{0,0}\ar[r]\ar[dr] &{\bf E}^{1,0}\ar[r]\ar[dr] &{\bf E}^{2,0}\ar[r]\ar[dr] &{\bf E}^{3,0} \ar[r]\ar[dr] &{\bf E}^{4,0}\ar[r]\ar[dr] &{\bf E}^{5,0}\ar[r] &{\bf E}^{6,0}\\
&{\bf E}^{1,1} \ar[ur]\ar[r]\ar[dr] &{\bf E}^{2,1}\ar[r]\ar[dr] \ar[ur]            &  {\bf E}^{3,1}\ar[r]\ar[dr]\ar[ur]  & {\bf E}^{4,1} \ar[r]\ar[ur]   &{\bf E}^{5,1}\ar[ur] & \\
 &               & {\bf E}^{2,2} \ar[r] \ar[ur] \ar[dr] & {\bf E}^{3,2} \ar[ur] \ar[r] & {\bf E}^{4,2}\ar[ur]&&\\
&&& {\bf E}^{3,3}\ar[ur]&&&}
$$

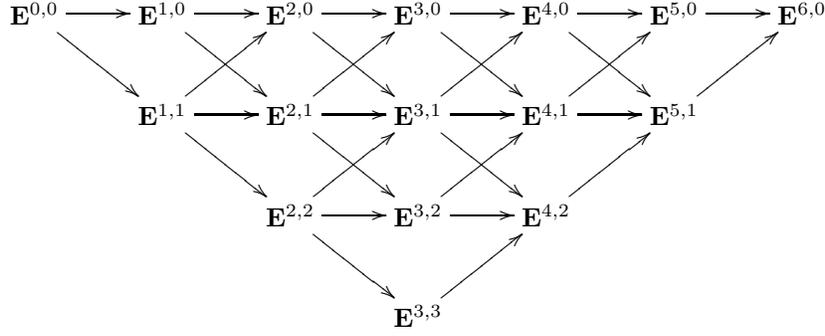
\centerline{Figure 1.}

Now, we shall introduce four operators which help us to describe the action of the symplectic Ricci curvature tensor field acting on symplectic spinor valued exterior differential forms. For $r=0,\ldots, 2l,$ $\alpha \otimes s \in \bigwedge^r
\mathbb{V}^*\otimes {\bf S}$ and $\sigma \in \odot^2 \mathbb{V}^*,$ we set
\begin{eqnarray*}
X&:& \bigwedge^{r}\mathbb{V}^*\otimes {\bf S} \to
\bigwedge ^{r+1}\mathbb{V}^*\otimes {\bf S}, \, \mbox{  }X(\alpha \otimes
s):=\sum_{i=1}^{2l}\epsilon^i\wedge \alpha \otimes e_i.s,\\
Y&:& \bigwedge^{r}\mathbb{V}^{*}\otimes {\bf S} \to \bigwedge^{r-1}\mathbb{V}^*\otimes {\bf S}, \, \mbox{  } Y(\alpha \otimes s):=\sum_{i,j=1}^{2l} \omega^{ij}\iota_{e_i}\alpha \otimes e_j.s, \\
\Sigma^{\sigma}&:&\bigwedge^{r}\mathbb{V}^*\otimes {\bf S}\to \bigwedge^{r+1} \mathbb{V}^* \otimes {\bf S}, \, \mbox{  } \Sigma^{\sigma} (\alpha \otimes s):=\sum_{i,j=1}^{2l}{\sigma^i}_j\epsilon^j\wedge \alpha \otimes e_i.s \, \mbox{ and  }\\
\Theta^{\sigma}&:&\bigwedge^{r}\mathbb{V}^*\otimes {\bf S} \to \bigwedge^{r}\mathbb{V}^*\otimes {\bf S}, \, \mbox{   } \Theta^{\sigma} (\alpha \otimes s):=\sum_{i,j=1}^{2l} \alpha \otimes \sigma^{ij}e_{ij}.s \, \mbox{   }
\end{eqnarray*}
and extend it linearly. Here $\sigma_{ij}:=\sigma(e_i,e_j),$ $i,j =1,\ldots, 2l,$ and the contraction of an exterior form $\alpha \in \bigwedge^{\bullet}\mathbb{V}^*$ by a vector $v\in \mathbb{V}$ is denoted by $\iota_{v}\alpha.$  

{\bf Remark:} 
\begin{itemize}
\item[1)] One easily finds out that the operators are independent of the choice of a symplectic basis $\{e_i\}_{i=1}^{2l}.$
The operators $X$ and $Y$ are used to prove the Howe correspondence for $Mp(\mathbb{V},\omega)$ acting on 
$\bigwedge^{\bullet}\mathbb{V}^*\otimes {\bf S}$ via the representation $\rho.$ See Kr\'ysl \cite{KryslJRT} for details. 
\item[2)] The symmetric tensor $\sigma$ is an infinitesimal version of a part of the curvature of a Fedosov connection. This part is called symplectic Ricci curvature tensor field and will be introduced bellow. The operators $\Sigma^{\sigma}$ and $\Theta^{\sigma}$ will help us to describe the action of the symplectic Ricci curvature tensor field acting on symplectic spinor valued exterior differential forms.
\end{itemize}

In what follows, we shall write $\iota_{e_{ij}}\alpha$ instead of $\iota_{e_i} \iota_{e_j}\alpha,$ $i,j=1,\ldots, 2l,$ and similarly for higher number of contracting elements.
  
Using the Lemma 1, it is easy to compute that
\begin{eqnarray}
X^2(\alpha\otimes s)=-\frac{\imath}{2}\omega_{ij}\epsilon^i\wedge \epsilon^j \wedge \alpha \otimes s \mbox{  and }\qquad  Y^2(\alpha \otimes s)=\frac{\imath}{2}\omega^{ij}\iota_{e_{ij}}\alpha \otimes s
\end{eqnarray}
for any element $\alpha \otimes s \in \bigwedge^{\bullet}\mathbb{V}^*\otimes {\bf S}.$
 
In order to be able to use the operators $X$ and $Y$ in a geometric setting and some further reasons, we shall need the following 

{\bf Lemma 3:}
\begin{itemize}
\item[1)] The operators $X,$ $Y$ are $\tilde{G}$-equivariant wr. to the representation $\rho$ of $\tilde{G}.$ 
\item[2)]For $(i,j) \in \Xi_-,$ the operator $X$ is an isomorphism if restricted to ${\bf E}^{ij}.$\\
          For $(i,j) \in \Xi_+,$ the operator $Y$ is an isomorphism if restricted to ${\bf E}^{ij}.$
\end{itemize}
{\it Proof.} For the $\tilde{G}$-equivariance of $X$ and $Y,$ see Kr\'ysl \cite{KryslRarita}.
The fact that the mentioned restrictions are isomorphisms is proved in Kr\'ysl \cite{KryslJRT}.
$\Box$

In the next lemma, four relations are proved which will be used later in order to determine a superset of the image of a restriction of the symplectic Ricci curvature tensor field acting on symplectic spinor valued exterior differential forms. Often, we shall write $\Sigma$ and $\Theta$ simply instead of the more explicit $\Sigma^{\sigma}$ and $\Theta^{\sigma}.$ The symmetric tensor $\sigma$ is assumed to be chosen. The symbol $\{,\}$ denotes the anticommutator on $\mbox{End}(\bigwedge^{\bullet}\mathbb{V}^*\otimes {\bf S}).$ 

{\bf Lemma 4:} The following relations
\begin{eqnarray}
 \{\Sigma, X \}         &=& 0, \\
 \left[ \{  \Sigma,  Y\} , Y^2 \right]   &=& 0,\\
  \left[X,\Theta \right]&=& 2\imath\Sigma \, \mbox{ and } \\
 \left[\Theta, Y^2\right]&=&0
\end{eqnarray} 
hold on $\bigwedge^{\bullet}\mathbb{V}^*\otimes {\bf S}.$

{\it Proof.}
We shall prove these identities for $\alpha\otimes s \in \bigwedge^{i}\mathbb{V}^{*}\otimes {\bf S},$ $i=0,\ldots, 2l$ only. The statement then follows by linearity of the considered operators.
\begin{itemize}
\item[1)] Let us compute
\begin{eqnarray*}
(X\Sigma + \Sigma X)(\alpha \otimes s)&=&X({\sigma^i}_j\epsilon^j\wedge \alpha \otimes e_i\phi)+\Sigma(\epsilon^i\wedge \alpha \otimes e_i.s)\\
                                         &=& {\sigma^i}_j \epsilon^k \wedge \epsilon^j \wedge \alpha \otimes e_{ki}.s
                                         +{\sigma^j}_k\epsilon^k \wedge \epsilon^i \wedge \alpha \otimes e_{ji}.s\\
                                         &=&{\sigma^i}_k\epsilon^j \wedge \epsilon^k \wedge \alpha \otimes e_{ji}.s
                                          +{\sigma^i}_k\epsilon^k\wedge \epsilon^j \wedge \alpha \otimes e_{ij}.s\\
                                         &=&{\sigma^i}_k\epsilon^j\wedge \epsilon^k \wedge \alpha \otimes (e_{ji}-e_{ij}).s\\
                                         &=&-\imath {\sigma^i}_k \omega_{ji} \epsilon^j \wedge \epsilon^k \wedge \alpha \otimes s\\
                                         &=&\imath \sigma_{jk} \epsilon^j\wedge \epsilon^k \wedge \alpha \otimes s\\
                                         &=&0,
\end{eqnarray*}
where we have renumbered indices, used the Lemma 1 and the fact that $\sigma$ is symmetric. In what follows, we shall use similar procedures without mentioning it explicitly.
\item[2)] Let us compute 
\begin{multline*}
P(\alpha \otimes s):=\{\Sigma, Y\}(\alpha \otimes s)\\
	\begin{aligned}
                      &= Y({\sigma^i}_j\epsilon^j \wedge \alpha \otimes e_i.s)+\Sigma (\omega^{ij}\iota_{e_i}\alpha\otimes 																									 e_j.s)\\
															 &= {\sigma^i}_j\omega^{kl}\iota_{e_k}(\epsilon^j \wedge \alpha) \otimes e_{li}.s +  																                                      \omega^{ij}{\sigma^k}_l\epsilon^l\wedge \iota_{e_{i}} \alpha\otimes e_{kj}.s \\
															 &={\sigma^i}_j \omega^{kl}(\delta^j_k\alpha - \epsilon^j \wedge \iota_{e_k}\alpha)\otimes e_{li}.s+
															    \omega^{ij}{\sigma^k}_l \epsilon^l \wedge \iota_{e_{i}}\alpha \otimes e_{kj}.s\\	
															 &=\sigma^{il}\alpha\otimes e_{li}.s-{\sigma^i}_j\omega^{kl}\epsilon^j \wedge \iota_{e_k}\alpha\otimes 																								e_{li}.s+ \omega^{ij}{\sigma^k}_l \epsilon^l \wedge \iota_{e_{i}}\alpha \otimes e_{kj}.s\\
															 &=\sigma^{il}\alpha\otimes e_{li}.s-{\sigma^k}_l\omega^{ij}\epsilon^l\wedge \iota_{e_i}\alpha \otimes e_{jk}.s+
																\omega^{ij}{\sigma^k}_l \epsilon^l \wedge \iota_{e_{i}}\alpha \otimes e_{kj}.s\\
															 &=\sigma^{il}\alpha\otimes e_{li}.s	-\imath \omega^{ij}\omega_{kj}{\sigma^k}_l\epsilon^l\wedge \iota_{e_{i}}\alpha 																		\otimes s\\
															 &=\sigma^{il}\alpha\otimes e_{li}.s - \imath {\sigma^i}_j\epsilon^j \wedge \iota_{e_i}\alpha \otimes s.					
							\end{aligned}
\end{multline*}
Now, we use the derived prescription for $P$ and the equation (1) to compute
\begin{multline*}
 \left[P,2\imath Y^2 \right](\alpha\otimes s)= 2\imath PY^2(\alpha\otimes s) - 2\imath Y^2 P(\alpha \otimes s)\\
 \quad\begin{aligned}
													 =& -P(\omega^{ij}\iota_{e_{ij}}\alpha \otimes s)	 - 2\imath Y^2 (\sigma^{ij}\alpha\otimes e_{ji}.s - \imath 																								{\sigma^i}_j\epsilon^j\wedge \iota_{e_i}\alpha \otimes s)\\
													 =& -\omega^{ij}\sigma^{kl}\iota_{e_{ij}}\alpha \otimes e_{lk}.s + \imath \omega^{ij}{\sigma^k}_l\epsilon^l \wedge 																					\iota_{e_{kij}}	\alpha \otimes s\\
													 &+\sigma^{ij}\omega^{kl}\iota_{e_{kl}}\alpha \otimes e_{ij}.s - \imath {\sigma^i}_j 																																	\omega^{kl}\iota_{e_{kl}}(\epsilon^j \wedge  \iota_{e_i}\alpha )\otimes s\\
													 =& -\omega^{ij}\sigma^{kl}\iota_{e_{ij}}\alpha \otimes e_{kl}.s + \imath                                                   \omega^{ij}{\sigma^k}_l \epsilon^l \wedge 	\iota_{e_{kij}}	\alpha \otimes s\\
													 &+ \sigma^{ij}\omega^{kl}\iota_{e_{kl}}\alpha \otimes e_{ji}.s -\imath 
													 	\omega^{kl}{\sigma^i}_j(\delta^j_l\iota_{e_{ki}}\alpha - \delta^j_k \iota_{e_{li}}\alpha                              + \epsilon^j \wedge \iota_{e_{kli}}\alpha)\otimes s\\
													=& -\omega^{ij}\sigma^{kl}\iota_{e_{ij}}\alpha \otimes e_{kl}.s + \imath                                                   \omega^{ij}{\sigma^k}_l \epsilon^l \wedge 	\iota_{e_{kij}}	\alpha \otimes s\\
													   &+\sigma^{kl}\omega^{ij}\iota_{e_{ij}}\alpha\otimes e_{kl}.s-\imath \omega^{ij}{\sigma^{k}}_l\epsilon^l\wedge \iota_{e_{ijk}}\alpha\otimes s\\	
													 =& \, 0.
													 \end{aligned}
\end{multline*}
\item[3)] Due to the definition of $\Theta,$ we have
\begin{eqnarray*}
\left[X, \Theta \right](\alpha\otimes s) &=& \epsilon^k \wedge \alpha \otimes \sigma^{ij}e_{kij}.s- \epsilon^i \wedge \alpha \otimes \sigma^{jk}e_{jki}.s\\
&=&\epsilon^k\wedge \alpha \otimes \sigma^{ij}e_{kij}.s - \epsilon^k \wedge \alpha \otimes \sigma^{ij}e_{ijk}.s\\
&=& \sigma^{ij}\epsilon^k\wedge \alpha\otimes (e_{ikj}.s- \imath \omega_{ki}e_j.s - e_{ijk}.s)\\
&=& \sigma^{ij}\epsilon^k\wedge \alpha\otimes (e_{ijk}.s- \imath \omega_{kj}e_i.s - \imath \omega_{ki}e_j.s - e_{ijk}.s)\\
&=&  2 \imath \Sigma (\alpha \otimes s).
\end{eqnarray*}

\item[4)] This  relation follows easily from the definition of $\Theta$ and the relation (1).
$\Box$
\end{itemize}

In the next proposition, a superset of the image of $\Sigma$  and $\Theta$ restricted to ${\bf E}^{ij},$ for $(i,j) \in \Xi,$ is determined.

{\bf Proposition 5:} For $(i,j) \in \Xi,$ we have  
\begin{eqnarray*}
\Sigma_{|{\bf E}^{ij}}&:& {\bf E}^{ij} \to {\bf E}^{i+1,j-1} \oplus {\bf E}^{i+1,j} \oplus {\bf E}^{i+1,j+1}\, \mbox{ and  }\\
\Theta_{|{\bf E}^{ij}}&:& {\bf E}^{ij} \to  {\bf E}^{i+1,j-1} \oplus {\bf E}^{i+1,j} \oplus {\bf E}^{i+1,j+1}.
\end{eqnarray*}

{\it Proof.} 
\begin{itemize}
\item[1)]For $i=0,\ldots, l,$ let us choose an element $\psi =\alpha \otimes s \in {\bf E}^{ii}.$ Using the relation (3), we have
$0=[P,Y^2]\psi = (P Y^2 - Y^2 P)\psi = (\Sigma Y^3 + Y\Sigma Y^2 -Y^2\Sigma Y + Y^3\Sigma)\psi.$ Because $Y$ is $\tilde{G}$-equivariant (Lemma 3 item 1), decreasing the form degree of $\psi$ by one and there is no summand isomorphic to ${\bf E}^{ii}_+$ or ${\bf E}^{ii}_-$ in $\bigwedge^{i-1}\mathbb{V}^*\otimes {\bf S}$ (Remark bellow the Theorem 2), $Y\psi=0.$ Using this equation, we see that the first three summands in the above expression for $[P,Y^2]$ are zero. Therefore we have $0=Y^3\Sigma \psi.$  Because $Y$ is injective on ${\bf E}^{ij}$ for $(i,j)\in \Xi_+$ (Lemma 3 item 2), we see that $\Sigma \psi \in {\bf E}^{i+1,i-1} \oplus {\bf E}^{i+1,i} \oplus {\bf E}^{i+1,i+1}.$ 

Now, let us consider a general $(i,j) \in \Xi$ and $\psi \in  {\bf E}^{ij}.$ Let us take an element $\psi' \in {\bf E}^{jj}$ such that $\psi=X^{(i-j)}\psi'.$ This element exists because according to Lemma 3 item 2, the operator $X$ is an isomorphism when restricted to ${\bf E}^{ij}$ for $(i,j)\in \Xi_-.$
Because of the relation (2), we have $\Sigma\psi=\Sigma X^{(i-j)}\psi'=\pm X^{(i-j)}\Sigma \psi'.$ From the previous item, we know that  $\Sigma \psi' \in {\bf E}^{j+1,j-1} \oplus {\bf E}^{j+1,j} \oplus {\bf E}^{j+1,j+1}.$ Because $X$ is $\tilde{G}$-equivariant (Lemma 3 item 1) and the only summands in $\bigwedge^{i+1}\mathbb{V}^*\otimes {\bf S}$ isomorphic to ${\bf E}^{j+1,j-1} \oplus {\bf E}^{j+1,j} \oplus {\bf E}^{j+1,j+1}$ are those described in the formulation of this proposition (see the Remark bellow the Theorem 2), the statement follows.

\item[2)] For $i=0,\ldots, l,$ let us consider an element $\psi = \alpha \otimes s \in {\bf E}^{ii}.$ Using the relation (5), we have $0=[\Theta,Y^2]\psi = \Theta Y^2 \psi +  Y^2 \Theta \psi.$ Using similar reasoning to that one in the first item, we get $Y\psi=0.$ Using the expression for $[\Theta, Y^2]$ above, we get $Y^2\Theta \psi = 0$ and consequently,  $\Theta \psi \in {\bf E}^{ii}\oplus {\bf E}^{i,i-1}.$ Now, let us suppose $\psi \in {\bf E}^{ij}$ for $(i,j) \in \Xi.$ There exists an element $\psi' \in {\bf E}^{jj}$ such that $\psi = X^{(i-j)}\psi'$  (Lemma 3 item 2).
Using the relations (4) and (2), we have $\Theta \psi = \Theta X^{(i-j)}\psi'= X^{(i-j)}\Theta  \psi'$ if $i-j$ is even and  $(X^{(i-j)}\Theta - 2\imath X^{(i-j-1)}\Sigma)\psi'$ if $i-j$ is odd. Using the fact $\Sigma_{|{\bf E}^{ij}}: {\bf E}^{ij} \to {\bf E}^{i+1,j-1} \oplus {\bf E}^{i+1,j} \oplus {\bf E}^{i+1,j+1},$ the statement follows by similar lines of reasoning as in the first item.
$\Box$
\end{itemize}
  
\section{Metaplectic structures and symplectic curvature tensors}  

After we have  finished the algebraic part of the paper, let us start describing the geometric structure we shall be investigating.
  We begin with a recollection of results of Vaisman in \cite{Vaisman} and of Gelfand, Retakh and Shubin in \cite{GSR}.  
Let $(M,\omega)$ be a symplectic manifold and 
$\nabla$ be a symplectic torsion-free affine connection. By symplectic and torsion-free, we  mean $\nabla \omega =0$ and $T(X,Y):=\nabla_XY-\nabla_YX-[X,Y]=0$ for all $X,Y \in \mathfrak{X}(M),$ respectively.
Such connections are usually  called Fedosov connections. In what follows, we shall call the triple $(M,\omega,\nabla)$ Fedosov manifolds.    

 To fix our notation, let us recall the classical definition of the curvature tensor $R^{\nabla}$ of the connection $\nabla,$ we shall be using here. Let 
$$R^{\nabla}(X,Y)Z:=\nabla_X\nabla_Y Z - \nabla_Y \nabla_X Z - \nabla_{[X,Y]}Z$$ for
$X,Y,Z \in \mathfrak{X}(M).$  

Let us choose a local symplectic frame $\{e_i\}_{i=1}^{2l}$ over an open subset $U\subseteq M.$
We shall often write expressions in which indices $i,j,k,l$ e.t.c. occur. We will implicitly mean $i,j,k, l$ are running from $1$ to $2l$ without mentioning it explicitly.
We set $$R_{ijkl}:=\omega(R(e_k,e_l)e_j,e_i).$$ Let us mention that we are using the convention of Vaisman \cite{Vaisman} which is different from that one used in Habarmann, Habermann \cite{HH}.

 From the symplectic curvature tensor field $R^{\nabla}$, we can build the symplectic Ricci curvature tensor field  $\sigma^{\nabla}$ defined by the 
classical formula

$$\sigma^{\nabla}(X,Y):=\mbox{Tr}(V \mapsto R^{\nabla}(V,X)Y)$$ for each $X,Y \in \mathfrak{X}(M)$ (the variable $V$ denotes a vector field on $M$). For the chosen frame and $i,j=1,\ldots, 2l$, we set
$$\sigma_{ij}:=\sigma^{\nabla}(e_i,e_j).$$

Further, let us define
\begin{eqnarray}
2(l+1)\widetilde{\sigma}^{\nabla}_{ijkl}&:=&\omega_{il}\sigma_{jk}-\omega_{ik}\sigma_{jl}+\omega_{jl}\sigma_{ik}-\omega_{jk}\sigma_{il}+2\sigma_{ij}\omega_{kl},\\
\widetilde{\sigma}^{\nabla}(X,Y,Z,V)&:=&\widetilde{\sigma}_{ijkl}X^iY^jZ^kV^l\, \mbox{ and} \nonumber\\
W^{\nabla}&:=&R^{\nabla}-\widetilde{\sigma}^{\nabla}
\end{eqnarray}
for  local vector fields $X=X^ie_i,$ $Y=Y^je_j,$ $Z=Z^ke_k$ and $V=V^le_l.$
We will call the tensor field $\tilde{\sigma}$ the extended symplectic Ricci curvature tensor field and $W^{\nabla}$ the symplectic Weyl curvature tensor field.
These tensor fields were already introduced in Vaisman \cite{Vaisman}. We shall often drop the index $\nabla$ in the previous expressions. Thus, we shall often write $W,$ $\sigma$ and $\widetilde{\sigma}$ instead of
$W^{\nabla},$ $\sigma^{\nabla}$ and $\widetilde{\sigma}^{\nabla},$ respectively.

In the next lemma, the symmetry of $\sigma$ is stated.

{\bf Lemma 6:} The symplectic Ricci curvature tensor field $\sigma$ is symmetric.

{\it Proof.}   See Vaisman \cite{Vaisman}. $\Box$




Let us start describing the geometric structure with help of which the action of the symplectic twistor operators  are defined.  This structure, called  metaplectic, is a precise 
symplectic analogue of the notion of a spin structure in  the Riemannian geometry.
For a symplectic manifold $(M^{2l}, \omega)$  of dimension $2l,$
let us denote the bundle of symplectic reperes in $TM$ by
$\mathcal{P}$ and  the foot-point projection of $\mathcal{P}$ onto
$M$ by $p.$ Thus $(p:\mathcal{P}\to M, G),$ where $G\simeq
Sp(2l,\mathbb{R}),$ is a principal $G$-bundle over $M$. As in
the subsection 2, let $\lambda: \tilde{G}\to G$ be a member
of the isomorphism class of the non-trivial two-fold coverings of
the symplectic group $G.$ In particular, $\tilde{G}\simeq
Mp(2l,\mathbb{R}).$ Further, let us consider a principal
$\tilde{G}$-bundle $(q:\mathcal{Q}\to M, \tilde{G})$ over the
symplectic manifold $(M,\omega).$ We call a pair
$(\mathcal{Q},\Lambda)$   metaplectic structure if  $\Lambda:
\mathcal{Q} \to \mathcal{P}$ is a surjective bundle homomorphism
over the identity on $M$ and if the following diagram,
$$\begin{xy}\xymatrix{
\mathcal{Q} \times \tilde{G} \ar[dd]^{\Lambda\times \lambda} \ar[r]&   \mathcal{Q} \ar[dd]^{\Lambda} \ar[dr]^{q} &\\
                                                            & &M\\
\mathcal{P} \times G \ar[r]   & \mathcal{P} \ar[ur]_{p} }\end{xy}$$
with the
horizontal arrows being respective actions of the displayed groups, commutes.
See, e.g.,  Habermann,  Habermann \cite{HH} and Kostant \cite{Kostant2} for
details on   metaplectic structures. Let us only remark, that typical examples of symplectic manifolds admitting a metaplectic structure are cotangent bundles of orientable manifolds (phase spaces), Calabi-Yau manifolds and complex projective spaces $\mathbb{CP}^{2k+1}$, $k \in \mathbb{N}_0.$

Let us denote the  vector bundle  associated to the introduced principal $\tilde{G}$-bundle
$(q:\mathcal{Q}\to M,\tilde{G})$ via the representation $meta$ on ${\bf S}$ by $\mathcal{S}.$ We shall call
this associated vector bundle  symplectic spinor bundle. Thus, we have $\mathcal{S}=\mathcal{Q}\times_{meta}{\bf S}.$ The sections $\phi \in \Gamma(M,\mathcal{S}),$ will be called  symplectic spinor fields.
Let us denote the space of symplectic valued exterior differential forms $\Gamma(M,\mathcal{Q}\times_{\rho}(\bigwedge^{\bullet}\mathbb{V}^*\otimes{\bf S}))$ by $\Omega^{\bullet}(M,\mathcal{S})$ and call it the space of  symplectic spinor valued forms simply. 
Further for $(i,j)\in \mathbb{Z}\times \mathbb{Z},$  we define the  associated vector bundles $\mathcal{E}^{ij}$ by the prescription $\mathcal{E}^{ij}:=\mathcal{Q}\times_{\rho} {\bf E}^{ij}.$   

Because the operators $X,Y$ are $\tilde{G}$-equivariant (Lemma 3 item 1), they lift to operators
acting on sections of the corresponding associated vector bundles.
We shall use the same symbols as for the defined operators as for
their "lifts" to the associated vector bundle structure. Because for each $i=0,\ldots 2l,$ the decomposition 
$\bigwedge^i\mathbb{V}^* \otimes {\bf S} \simeq \bigoplus_{j,(i,j)\in \Xi} {\bf E}^{ij}$ is multiplicity free (see the Remark bellow the Theorem 2), there exist uniquely defined projections $p^{ij}:\Omega^i(M,\mathcal{S}) \to \Gamma(M,\mathcal{E}^{ij}),$ $(i,j)\in \mathbb{Z}\times \mathbb{Z}.$

Now, let us suppose that   $(M,\omega)$ is equipped with a Fedosov connection $\nabla$. The   connection $\nabla$ determines the associated principal bundle connection $Z$
on the principal bundle $(p:\mathcal{P}\to M, G).$ 
This connection lifts to a principal bundle connection on  the principal bundle
$(q:\mathcal{Q}\to M, \tilde{G})$ and defines the associated covariant derivative on the symplectic bundle $\mathcal{S},$ which we shall denote by $\nabla^S$ and call it the  symplectic spinor covariant derivative. See Habermann, Habermann \cite{HH} for details. The symplectic spinor covariant derivative induces the exterior symplectic spinor derivative $d^{\nabla^S}$ acting on $\Omega^{\bullet}(M,\mathcal{S}).$
The curvature tensor field $R^{\Omega^{\bullet}(M,\mathcal{S})}$ acting on the symplectic spinor valued forms is given by the classical formula
$$R^{\Omega^{\bullet}(M,\mathcal{S})}:=d^{\nabla^S} d^{\nabla^S}.$$  

In the next theorem, a superset of the image of $d^{\nabla^S}$ restricted to $\Gamma(M,\mathcal{E}^{ij}),$ $(i,j) \in \Xi,$ is determined. 

{\bf Theorem  7:} Let $(M,\omega,\nabla)$ be a Fedosov manifold admitting a metaplectic structure. Then for the exterior symplectic spinor derivative $d^{\nabla^S},$ we have
$$d^{\nabla^S}_{|\Gamma(M,\mathcal{E}^{ij})}: \Gamma(M,\mathcal{E}^{ij})\to \Gamma(M,\mathcal{E}^{i+1,j-1}\oplus \mathcal{E}^{i+1,j}\oplus \mathcal{E}^{i+1,j+1}),$$ where $(i,j) \in \Xi.$

{\it Proof.} See Kr\'ysl \cite{KryslSVF}. $\Box$

{\bf Remark:}
From the proof of the theorem, it is easy to see that it can be extended to the case $(M,\omega)$ is presymplectic and the symplectic connection $\nabla$ has a non-zero torsion.  
For $l=3$ and any $(i,j)\in \Xi,$ the mappings $d^{\nabla^S}$ restricted to $\Gamma(M,\mathcal{E}^{ij})$ are displayed as arrows at the Figure 1 above. (The exterior covariant derivative $d^{\nabla^S}$ maps $\Gamma(M,\mathcal{E}^{ij})$ into the three "neighbor" subspaces.)  

\subsection{Curvature tensor on symplectic spinor valued forms and the complex of symplectic twistor operators}

Let $(M,\omega, \nabla)$ be a Fedosov manifold admitting a metaplectic structure $(\mathcal{Q},\Lambda).$ 
In the next lemma, the action of $R^S:=d^{\nabla^S}\circ \nabla^S$ on the space of symplectic spinors  fields  is described using just the symplectic curvature tensor field $R$ of $\nabla.$
 
{\bf Lemma 8:} Let $(M,\omega,\nabla)$ be a 
Fedosov manifold  admitting a metaplectic structure.
Then for a symplectic spinor field $\phi \in \Gamma(M,\mathcal{S}),$ we have
$$R^S\phi=\frac{\imath}{2} {R^{ij}}_{kl}\epsilon^k\wedge \epsilon^l \otimes e_i.e_j.\phi.$$
  
{\it Proof.} See  Habermann, Habermann \cite{HH} pp. 42. $\Box.$
 
 
For our convenience, let us set $m_i:=i$ for $i=0,\ldots, l$ and $m_i:=2l-i$ for $i=l+1,\ldots, 2l.$
Now, we can define the symplectic twistor operators, which we shall need to introduce the mentioned complex.
For $i=0,\ldots, 2l-1,$ we set
$$T_i:\Gamma(M,\mathcal{E}^{im_i})\to \Gamma(M,\mathcal{E}^{i+1,m_{i+1}}), \quad \mbox{  } T_i:=p^{i+1, m_{i+1}}d^{\nabla^S}_{|\Gamma(M,\mathcal{E}^{im_i})}$$ and call these operators  {\it symplectic twistor operators.} 
Informally, one can say that the operators are going on the edge of the triangle at the Figure 1. 
Let us notice that $Y(\nabla^S-T_0)$ is, up to a nonzero scalar multiple, the so called symplectic Dirac operator introduced by K. Habermann in \cite{KH}.

{\bf Theorem 9:} Let $(M^{2l},\omega,\nabla)$ be a Fedosov manifold admitting a metaplectic structure. If $l\geq 2$ and the symplectic Weyl tensor field $W^{\nabla}=0,$  then
$$0  \longrightarrow
  \Gamma(M,\mathcal{E}^{00}) \overset{T_0}{\longrightarrow} 
  \Gamma(M,\mathcal{E}^{11}) \overset{T_{1}}{\longrightarrow}  
  \cdots  \overset{T_{l-1}}{\longrightarrow}
  \Gamma(M,\mathcal{E}^{ll}) \longrightarrow 0 \mbox{   and}  
$$ 
$$0  \longrightarrow
  \Gamma(M,\mathcal{E}^{ll}) \overset{T_l}{\longrightarrow} 
  \Gamma(M,\mathcal{E}^{l+1,l+1}) \overset{T_{l+1}}{\longrightarrow}  
  \cdots  \overset{T_{2l-1}}{\longrightarrow}
  \Gamma(M,\mathcal{E}^{2l,2l}) \longrightarrow 0  
$$ 
are complexes.

{\it Proof.}
\begin{itemize}
\item[1)] In this item, we prove that for an element $\psi \in \Omega^{\bullet}(M,\mathcal{S}),$ 
$$ R^{\Omega^{\bullet}(M,\mathcal{S})} \psi = \frac{\imath}{l+1}(\imath X^2\Theta^{\sigma}  - X\Sigma^{\sigma})\psi.$$ 
For $\psi=\alpha \otimes \phi \in \Omega^{\bullet}(M,\mathcal{S}),$ we can write
\begin{multline*}
								 R^{\Omega^{\bullet}(M,\mathcal{S})}(\alpha \otimes \phi)     = d^{\nabla^{S}}d^{\nabla^{S}}(\alpha \otimes \phi)  = d^{\nabla^{S}}(d\alpha \otimes \phi + (-1)^{deg(\alpha)}\alpha \wedge \nabla^S \phi)\\
								 \begin{aligned}
                                              &= d^2\alpha \otimes \phi + (-1)^{deg(\alpha)+1}d\alpha \wedge \nabla^S \phi+(-1)^{deg(\alpha)}d\alpha\wedge\nabla^S \psi+\\
                                               &\, (-1)^{deg(\alpha)}(-1)^{deg(\alpha)}\alpha \wedge d^{\nabla^S}  \nabla^S\phi = \alpha \wedge \frac{\imath}{2} {R^{ij}}_{kl}\epsilon^k\wedge \epsilon^l \otimes e_{ij}.\phi\\
                                              &= \frac{\imath}{2} {R^{ij}}_{kl}\epsilon^k \wedge \epsilon^l \wedge \alpha \otimes e_{ij}.\phi,
\end{aligned}
\end{multline*} 
where we have used the Lemma 8.                    
Using this computation, the definition of the symplectic Weyl curvature tensor field $W^{\nabla}$ (Eqn. (7)), the definition of the extended symplectic Ricci curvature tensor field $\widetilde{\sigma}$ (Eqn. (6))  and the assumption $W^{\nabla}=0$, we get 
\begin{multline*}
-4(l+1)\imath R^{\Omega^{\bullet}(M,\mathcal{S})}(\alpha\otimes\phi)=2(l+1){R^{ij}}_{kl}\epsilon^k\wedge\epsilon^l\wedge \alpha \otimes e_{ij}.\phi\\
\begin{aligned}
&=2(l+1)({W^{ij}}_{kl}+{\mbox{$\widetilde{\sigma}^{ij}$}}_{kl})\epsilon^k\wedge\epsilon^l\wedge\alpha\otimes e_{ij}.s\\
&=2(l+1){\mbox{$\widetilde{\sigma}^{ij}$}}_{kl}\epsilon^k\wedge\epsilon^l\wedge \alpha \otimes e_{ij}.s\\
&=({\omega^i}_l{\sigma^{j}}_k-{\omega^{i}}_{k}{\sigma^j}_l+{\omega^{j}}_l{\sigma^{i}}_{k}-{\omega^{j}}_k{\sigma^{i}}_{l}+2\sigma^{ij}\omega_{kl})\epsilon^k\wedge\epsilon^l\wedge\alpha \otimes e_{ij}.\phi\\
&=(4{\omega^i}_l{\sigma^j}_k\epsilon^k\wedge\epsilon^l\wedge \alpha\otimes e_{ij}.\phi+2\sigma^{ij}\omega_{kl})\epsilon^k\wedge \epsilon^l \wedge \alpha \otimes e_{ij}.\phi\\
&=4\imath X^2(\alpha \otimes \sigma^{ij} e_{ij}.\phi)-4 X ({\sigma^j}_k\epsilon^k \wedge \alpha \otimes e_j.\phi)=(4\imath X^2\Theta^{\sigma}-4X\Sigma^{\sigma})\psi, 
\end{aligned}
\end{multline*}where we have used 
the relation (1) in the second last step. Extending the result by linearity, we get the statement of this item for arbitrary $\psi \in \Omega^{\bullet}(M,\mathcal{S}).$

\item[2)]Using the derived formula for $R^{\Omega^{\bullet}(M,\mathcal{S})},$ the Proposition 5, the $\tilde{G}$-equi\-variance of $X$ (Lemma 3 item 1) and the decomposition structure of $\bigwedge^{\bullet}\mathbb{V}^*\otimes {\bf S}$ (see the Remark bellow the Theorem 2), we see that for $(i,j)\in \Xi$ and an element $\psi \in \Gamma(M,\mathcal{E}^{ij}),$ the section $R^{\Omega^{\bullet}(M,\mathcal{S})} \psi \in \Gamma(M,\mathcal{E}^{i+2,j-1} \oplus \mathcal{E}^{i+2,j} \oplus \mathcal{E}^{i+2,j+1}).$ Thus especially, $p^{i+2,m_{i+2}}R^{\Omega^{\bullet}(M,\mathcal{S})} \psi = 0$  for \linebreak $i=0,\ldots, l-2,l, \ldots, 2l-2$ and $\psi \in \Gamma(M,\mathcal{E}^{im_i}).$  
For $i=0, \ldots, l-2,$ we get
\begin{multline*}
\begin{aligned}
0&=p^{i+2,i+2}R^{\Omega^{\bullet}(M, \mathcal{S})}=p^{i+2,i+2}d^{\nabla^S}d^{\nabla^S}\\
&=p^{i+2,i+2}d^{\nabla^S}(p^{i+1,0} + \ldots + p^{i+1,i+1})d^{\nabla^S}\\
&=p^{i+2,i+2}d^{\nabla^S}p^{i+1,0}d^{\nabla^S}+\ldots + p^{i+2,i+2}d^{\nabla^S}p^{i+1,i+1}d^{\nabla^S}\\
&= T_{i+1}T_{i},
\end{aligned}
\end{multline*}
where we have used the Theorem 7 in the   last step.
Similarly, one proceeds in the case $i=l,\ldots, 2l-2.$ 
\end{itemize}
$\Box$

{\bf Corollary 10.} Let $(M,\omega, \nabla)$ be a Fedosov manifold admitting a metaplectic structure. If $l \geq 2$ and the symplectic Weyl tensor field $W^{\nabla}=0,$ then
$$ 0\longrightarrow \Gamma(M,\mathcal{E}^{00}) \overset{T_0}{\longrightarrow}\cdots \overset{T_{l-2}}{\longrightarrow} \Gamma(M,\mathcal{E}^{l-1,l-1})\overset{T_lT_{l-1}}{\longrightarrow}$$ 
$$\overset{T_lT_{l-1}}{\longrightarrow} \Gamma(M,\mathcal{E}^{l+1,l+1}) \overset{T_{l+1}}{\longrightarrow} \cdots \overset{T_{2l-1}}{\longrightarrow} \Gamma(M,\mathcal{E}^{2l,2l}) \longrightarrow 0
$$
is a complex.

{\it Proof.} Follows easily from the Theorem 9. $\Box$
 
The question of the existence of a symplectic connection with vanishing symplectic Weyl curvature tensor field was treated, e.g., in Cahen, Gutt, Rawnsley \cite{Cahen}. These connections are called connections of Ricci type.
For instance it is known that if a compact simply connected  symplectic manifold  $(M,\omega)$ admits a connection of Ricci type, then $(M,\omega)$ is affinely symplectomorphic to a $\mathbb{P}^n \mathbb{C}$ with the symplectic form, given by the standard complex structure and the Fubini-Study metric, and the Levi-Civita connection of this metric. 
Let us  refer an interested reader to the paper of Cahen, Gutt, Schwachh\"ofer \cite{CGS}, where also a relation of symplectic connections to contact projective geometries is treated.

Further research could be devoted to the investigation and the interpretation of the cohomology of the introduced complex and to the investigation of analytic properties of the introduced symplectic twistor operators.

\end{document}